\documentclass[a4paper,10pt]{elsart}
\usepackage{amsmath,amsfonts,amssymb,latexsym}
\newcommand{\fd}{\mathbb{F}}

\newcommand{\Z}{\mathbb{Z}}

\newtheorem{theorem}{Theorem}
\newtheorem{lemma}[theorem]{Lemma}
\newtheorem{corollary}[theorem]{Corollary}
\newtheorem{proposition}[theorem]{Proposition}
\newtheorem{remark}[theorem]{Remark}
\newcommand{\dis}{\displaystyle}

\begin{document}

\begin{frontmatter}

\title{Simple Recursive Formulas Generating Power Moments of Kloosterman Sums}

\author{Dae San Kim}

\address{Department of Mathematics, Sogang University, Seoul 121-742, South Korea}

\begin{abstract}
In this paper, we construct four binary linear codes closely
connected with  certain  exponential sums  over the finite field
$\fd_q$  and $\fd_q - \{0,1\}$. Here $q$ is a power of two. Then
we obtain four recursive formulas for the power moments of
Kloosterman sums in terms of the frequencies of weights in the
codes. This is done via Pless power moment identity and by
utilizing the explicit expressions of the exponential sums
obtained earlier.
\end{abstract}

\thanks{This work was supported by grant No. R01-2006-000-11176-0 from the
Basic Research Program of the Korea Science and Engineering
Foundation.}
\thanks{Email address : dskim@sogang.ac.kr.}


\end{frontmatter}

\section{Introduction and Notations}
Let  $\psi$ be a nontrivial additive character of the finite field
$\fd_q$ with $q = p^r$ elements ($p$ a prime). Then the Kloosterman
sum $K(\psi;a)$ is defined as
\[
K(\psi;a) = \sum_{\alpha \in \fd_q^*} \psi(\alpha + a\alpha^{-1})
(a \in \fd_q^*).
\]
The Kloosterman sum was introduced in 1926 to give an estimate for
the Fourier coefficients of modular forms. It has also been studied
to solve  various problems in coding theory and cryptography over
finite fields of characteristic two.

For each nonnegative integer $h$, by $MK(\psi)^h$  we will denote
the $h$-th moment of the Kloosterman sum $K(\psi;a)$. Namely, it
is given by
\[
MK(\psi)^h = \sum_{a \in \fd_q^*} K(\psi;a)^h.
\]
If $\psi = \lambda$ is the canonical additive character of $\fd_q$,
then $MK(\lambda)^h$ will be simply denoted by $MK^h$.

From now on, let us assume that $q=2^r$.   Carlitz \cite{L1}
evaluated $MK^{h}$, for $h\leq 4$. Recently, Moisio was able to find
explicit expressions of $MK^h$, for the other values of $h$ with $h
\leq 10$ (cf.\cite{M1}). This was done, via Pless power moment
identity, by connecting moments of Kloosterman sums and the
frequencies of weights in the binary Zetterberg code of length
$q+1$, which were known by the work of Schoof and Vlugt in
\cite{RS}.

In this paper, we will produce four recursive formulas generating
power moments of Kloosterman sums. To do that, we construct four
binary linear codes  closely connected with certain exponential
sums(cf. Theorem 2). It is amusing that those exponential sums have
been used repeatedly in our derivation of recursive formulas
generating power moments of Kloosterman sums (cf. \cite{Kim2},
\cite{Kim3}).

Theorem 1 of the following(cf. (1), (2), (4)-(9)) is the main result
of this paper. Henceforth, we agree that the binomial coefficient
${b \choose a}=0$, if $a>b$ or $a<0$.

\begin{theorem} Let $q=2^r$. Then we have the following.\\
(a)For $r \geq 3$, and $h=1,2,\ldots,$
\begin{align}
\begin{split}
MK^h =& \sum_{l=0}^{h-1} (-1)^{h+l+1} {h \choose l} (q-1)^{h-l}
MK^l \\
      &+ q \sum_{j=0}^{min \{q-2, ~h\}} (-1)^{h+j} C_{1,j}
      \sum_{t=j}^h t! S(h,t) 2^{h-t} {q-2-j \choose q-2-t},
\end{split}
\end{align}
where $\{C_{1,j}\}_{j=0}^{q-2}$ is the weight distribution of the
 binary linear code $C_1$ given by
\begin{equation}
C_{1,j} = \sum \prod_{tr(\beta^{-1})=0} {2 \choose \nu_\beta}
(j=0,\ldots,q-2).
\end{equation}
Here the sum is over all the sets of nonnegative integers
$\{\nu_\beta\}_{tr(\beta^{-1})=0}$ satisfying
$\dis\sum_{tr(\beta^{-1})=0} \nu_\beta = j$ and
$\dis\sum_{tr(\beta^{-1})=0} \nu_\beta \beta = 0$. In addition,
$S(h,t)$ is the Stirling number of the second kind defined by
\begin{equation}
S(h,t) = \frac{1}{t!}\sum_{j=0}^t (-1)^{t-j} {t \choose j} j^h.
\end{equation}\\
(b) For $r \geq 3$, and $h=1,2,\ldots,$
\begin{align}
\begin{split}
MK^h =& \sum_{l=0}^{h-1} (-1)^{h+l+1} {h \choose l} (q-1)^{h-l}
MK^l \\
      &+ q \sum_{j=0}^{min \{\frac{q}{2}-1, ~h\}} (-1)^{h+j} C_{2,j}
      \sum_{t=j}^h t! S(h,t) 2^{2h-t} {\frac{q}{2}-1-j \choose \frac{q}{2}-1-t},
\end{split}
\end{align}
where $\{C_{2,j}\}_{j=0}^{\frac{q}{2}-1}$ is the weight
distribution of the  binary linear code $C_2$, and
\begin{align}
\begin{split}
&C_{2,j}(0 \leq j \leq \frac{q}{2}-1) ~~ \text{is equal to the
number of all the sets of nonnegative integers}\\
&\{\nu_\beta\}_{tr(\beta^{-1})=0} ~~\text{satisfying} ~~\nu_\beta
= 0 ~~or ~~1, \dis\sum_{tr(\beta^{-1})=0} \nu_\beta = j ~~and
\dis\sum_{tr(\beta^{-1})=0} \nu_\beta \beta = 0.
\end{split}
\end{align}
(c) For $h=1,2,\ldots,$
\begin{align}
\begin{split}
MK^h =& - \sum_{l=0}^{h-1} {h \choose l} (q+1)^{h-l}MK^l \\
      &+ q \sum_{j=0}^{min \{q, ~h\}} (-1)^j C_{3,j}
      \sum_{t=j}^h t! S(h,t) 2^{h-t} {q-j \choose q-t},
\end{split}
\end{align}
where $\{C_{3,j}\}_{j=0}^q$ is the weight distribution of the
 binary linear code $C_3$ given by
\begin{equation}
C_{3,j} = \sum \prod_{tr(\beta^{-1})=1} {2 \choose \nu_\beta}
(j=0,\ldots,q).
\end{equation}
Here the sum is over all the sets of nonnegative integers
$\{\nu_\beta\}_{tr(\beta^{-1})=1}$ satisfying
$\dis\sum_{tr(\beta^{-1})=1} \nu_\beta = j$ and
$\dis\sum_{tr(\beta^{-1})=1} \nu_\beta \beta = 0$.\\
(d) For $h=1,2,\ldots,$
\begin{align}
\begin{split}
MK^h =& - \sum_{l=0}^{h-1} {h \choose l} (q+1)^{h-l}MK^l \\
      &+ q \sum_{j=0}^{min \{\frac{q}{2}, ~h\}} (-1)^j C_{4,j}
      \sum_{t=j}^h t! S(h,t) 2^{2h-t} {\frac{q}{2}-j \choose \frac{q}{2}-t},
\end{split}
\end{align}
where $\{C_{4,j}\}_{j=0}^{\frac{q}{2}}$ is the weight distribution
of the
 binary linear code $C_4$, and
\begin{align}
\begin{split}
&C_{4,j}(0 \leq j \leq \frac{q}{2}) ~~ \text{is equal to the
number of all the sets of nonnegative integers}\\
&\{\nu_\beta\}_{tr(\beta^{-1})=1} ~~\text{satisfying} ~~\nu_\beta
= 0 ~~or ~~1, \dis\sum_{tr(\beta^{-1})=1} \nu_\beta = j ~~and
\dis\sum_{tr(\beta^{-1})=1} \nu_\beta \beta = 0.
\end{split}
\end{align}
\end{theorem}
Before we proceed further, we will fix the notations that will be
used throughout this paper:
\begin{itemize}
 \item [] $q = 2^r$ ($r \in \Z_{>0}$),\\
 \item [] $\fd_q$ = the finite field with $q$ elements,\\
 \item [] $tr(x)=x+x^2+\cdots+x^{2^{r-1}}$ the trace function $\fd_q
\rightarrow \fd_2$,\\
 \item [] $\lambda(x) = (-1)^{tr(x)}$ the canonical additive
character of $\fd_q$.
\end{itemize}

\section{Two exponential sums}
Let  $\fd_2^+, \fd_q^+$ denote the additive groups of the fields
$\fd_2,\fd_q$, respectively. Then, with  $\Theta(x) = x^2+x$
denoting the Artin-Schreier operator in characteristic two, we
have the following exact sequence of groups:
\begin{equation*}
0 \rightarrow \fd_2^+ \rightarrow \fd_q^+ \rightarrow
\Theta(\fd_q) \rightarrow 0.
\end{equation*}
Here the first map is the inclusion and the second is  given by $x
\mapsto \Theta(x) = x^2+x$.  So
\begin{equation}
\Theta(\fd_q) = \{\alpha^2 + \alpha \mid  \alpha \in \fd_q \},~
and ~~[\fd_q^+ : \Theta(\fd_q)] = 2.
\end{equation}

The next theorem is the key to the results of this paper.
\begin{theorem}[\cite{Kim1}]
 Let $\lambda$  be the canonical additive character of $\fd_q$, and let $a \in \fd_q^*$. Then
\begin{align*}
 &(a) \sum_{\alpha \in
 \fd_q-\{0,1\}}\lambda(\frac{a}{\alpha^2+\alpha})= K(\lambda;a)-1,~~~~~~~~~~~~~~~~~~~~~~~~~~~~~~~~~~~~~~~\\
 &(b)\sum_{\alpha \in
\fd_q}\lambda(\frac{a}{\alpha^2+\alpha+b}) = -K(\lambda;a)-1, ~if
~~x^2+x+b (b \in \fd_q) ~~\text{is
irreducible over} ~\fd_q,\\
& ~~\text{or equivalently if} ~~b \in \fd_q\setminus\Theta(\fd_q)
(cf.(10)).
\end{align*}
\end{theorem}

\section{Construction of codes}
Here we will construct four binary linear codes $C_1$ of length
$N_1=q-2(q \geq 4)$, $C_2$ of length $N_2=\frac{q}{2}-1(q \geq
4)$, $C_3$ of length $N_3=q$, and $C_4$ of length
$N_4=\frac{q}{2}$, which are closely connected with the
exponential sums (a), (b) in Theorem 2.

Let $b \in \fd_q\setminus \Theta(\fd_q)$ be fixed, and let
$\gamma_0=0, \gamma_1,\ldots,\gamma_{\frac{q}{2}-1}$ be a fixed
ordering of the elements in $\Theta(\fd_q)= \{\alpha^2+\alpha \mid
\alpha \in \fd_q\}$. Then
\begin{align*}
& \Theta(\fd_q) = \{
\gamma_0,\gamma_1,\ldots,\gamma_{\frac{q}{2}-1}\} = \{\beta \in
\fd_q \mid tr \beta = 0\},\\
& b + \Theta(\fd_q) = \{
b+\gamma_0,b+\gamma_1,\ldots,b+\gamma_{\frac{q}{2}-1}\} = \{\beta
\in \fd_q \mid tr\beta = 1\},\\
&\fd_q = \{\gamma_0,\gamma_1,\ldots,\gamma_{\frac{q}{2}-1},
b+\gamma_0,b+\gamma_1,\ldots,b+\gamma_{\frac{q}{2}-1}\}.
\end{align*}
Also, we put
\begin{align}
&v_1 =
(\frac{1}{\gamma_1},\frac{1}{\gamma_2},\ldots,\frac{1}{\gamma_{\frac{q}{2}-1}},\frac{1}{\gamma_1},\frac{1}{\gamma_2},\ldots,\frac{1}{\gamma_{\frac{q}{2}-1}})
\in \fd_q^{N_1},\\
&v_2 =
(\frac{1}{\gamma_1},\frac{1}{\gamma_2},\ldots,\frac{1}{\gamma_{\frac{q}{2}-1}})
\in \fd_q^{N_2},\\
&v_3 =
(\frac{1}{b+\gamma_0},\frac{1}{b+\gamma_1},\ldots,\frac{1}{b+\gamma_{\frac{q}{2}-1}},\frac{1}{b+\gamma_0},\frac{1}{b+\gamma_1},\ldots,\frac{1}{b+\gamma_{\frac{q}{2}-1}})
\in \fd_q^{N_3},\\
&v_4 =
(\frac{1}{b+\gamma_0},\frac{1}{b+\gamma_1},\ldots,\frac{1}{b+\gamma_{\frac{q}{2}-1}})
\in \fd_q^{N_4}.
\end{align}

\begin{remark}
For each $i=1,2,3,4$, and each $\beta \in \fd_q$, let $n_i(\beta)$
denote the number of components with those equal to $\beta$ in the
vector $v_i$. Then it is obvious that
\begin{align}
& n_1(\beta) = \left\{%
\begin{array}{ll}
    2, & \hbox{if $\beta \in \fd_q^*$, with ~$tr(\beta^{-1}) = 0$,} \\
    0, & \hbox{otherwise,} \\
\end{array}%
\right.\\
& n_2(\beta) = \left\{%
\begin{array}{ll}
    1, & \hbox{if $\beta \in \fd_q^*$, with ~$tr(\beta^{-1}) = 0$,} \\
    0, & \hbox{otherwise,} \\
\end{array}%
\right.\\
& n_3(\beta) = \left\{%
\begin{array}{ll}
    2, & \hbox{if ~$tr(\beta^{-1}) = 1$,} \\
    0, & \hbox{otherwise,} \\
\end{array}%
\right.\\
& n_4(\beta) = \left\{%
\begin{array}{ll}
    1, & \hbox{if ~$tr(\beta^{-1}) = 1$,} \\
    0, & \hbox{otherwise.} \\
\end{array}%
\right.
\end{align}
Then, for each  $i=1,2,3,4,$ the binary linear code $C_i$ is
defined as
\begin{equation}
C_i = \{u \in \fd_2^{N_i} \mid u \cdot v_i = 0\},
\end{equation}
where the dot denotes the usual inner product in $\fd_q^{N_i}$.
\end{remark}

The following Delsarte's theorem is well-known.

\begin{theorem}[\cite{FJ}]
Let  $B$ be a linear code over $\fd_q$.  Then
\[
(B|_{\fd_2})^\bot = tr(B^\bot).
\]
\end{theorem}

In view of this theorem, the duals  $C_i^\bot (i=1,2,3,4)$ are
given by
\begin{align}
&C_1^\bot = \{c_1(a) =
(tr\frac{a}{\gamma_1},tr\frac{a}{\gamma_2},\ldots,tr\frac{a}{\gamma_{\frac{q}{2}-1}},tr\frac{a}{\gamma_1},tr\frac{a}{\gamma_2},\ldots,tr\frac{a}{\gamma_{\frac{q}{2}-1}})
\mid a \in \fd_q \},\\
&C_2^\bot = \{c_2(a) =
(tr\frac{a}{\gamma_1},tr\frac{a}{\gamma_2},\ldots,tr\frac{a}{\gamma_{\frac{q}{2}-1}})
\mid a \in \fd_q\},\\
&C_3^\bot = \{c_3(a) =
(tr\frac{a}{b+\gamma_0},\ldots,tr\frac{a}{b+\gamma_{\frac{q}{2}-1}},tr\frac{a}{b+\gamma_0},\ldots,tr\frac{a}{b+\gamma_{\frac{q}{2}-1}})
\mid a \in \fd_q\},\\
&C_4^\bot = \{c_4(a) =
(tr\frac{a}{b+\gamma_0},tr\frac{a}{b+\gamma_1},\ldots,tr\frac{a}{b+\gamma_{\frac{q}{2}-1}})
\mid a \in \fd_q\}.
\end{align}

(a) and (b) can be respectively proved just as Theorem 11 in
\cite{Kim2} and Theorem 13 in \cite{Kim3}.

\begin{theorem}
(a) For $q=2^r$, with $r \geq 3$, the map $\fd_q \rightarrow
C_i^\bot$$(a \mapsto c_i(a))$, for $i=1,2$, is an $\fd_2$-linear
isomorphism.\\
(b) For every $q=2^r$, the map $\fd_q \rightarrow C_i^\bot$$(a
\mapsto c_i(a))$, for $i=3,4$, is an $\fd_2$-linear isomorphism.
\end{theorem}

\begin{remark}
Note that, for $i=1,2$, and $q=4$, the kernel of the map  $\fd_4
\rightarrow C_i^\bot$$(a \mapsto c_i(a))$ is  $\fd_2$.
\end{remark}

\section{Power moments of Kloosterman sums}
In this section, we will be able to find, via Pless power moment
identity, a recursive formula for the power moments of Kloosterman
sums in terms of the frequencies of weights in  $C_i$, for each
$i=1,2,3,4.$

\begin{theorem}[Pless power moment identity]
Let $B$ be an $q$-ary $[n,k]$ code, and let $B_i$ (resp.
$B_i^\bot$) denote the number of codewords of weight $i$ in
$B$(resp. in $B^\bot$).  Then, for $h=0,1,2,\ldots,$
\begin{equation}
\sum_{j=0}^n j^h B_j=\sum_{j=0}^{min\{n,h\}} (-1)^j B_j^\bot
 \sum_{t=j}^h t! S(h,t) q^{k-t}(q-1)^{t-j}{n-j\choose n-t},
\end{equation}
where $S(h,t)$ is the Stirling number of the second kind defined
in (4).
\end{theorem}

\begin{lemma}
For $i=1,2,3,4,$ and $a \in \fd_q^*$, the Hamming weight $w(c_i(a))$
(cf. (20)-(23)) can be expressed as follows:
\begin{align}
&(a) ~w(c_1(a)) = \frac{1}{2}(q-1-K(\lambda;a)),~~~~~~~~~~~~~~~~~~~~~~~~~~~~~~~~~~~~~~~~\\
&(b) ~w(c_2(a)) = \frac{1}{4}(q-1-K(\lambda;a)),\\
&(c) ~w(c_3(a)) = \frac{1}{2}(q+1+K(\lambda;a)),\\
 &(d) ~w(c_4(a)) =
\frac{1}{4}(q+1+K(\lambda;a)).
\end{align}
\end{lemma}
\begin{pf*}{Proof.}
We will just show (c) and (d),  as (a) and (b) can be proved in a
similar manner.
\begin{align*}
w(c_3(a)) &= \frac{1}{2} \cdot 2 \cdot \sum_{j=0}^{\frac{q}{2}-1}
(1-(-1)^{tr(\frac{a}{b+\gamma_j})}) (cf. (22))\\
&= \frac{1}{2} (q-2 \sum_{j=0}^{\frac{q}{2}-1}
\lambda(\frac{a}{b+\gamma_j}))\\
&= \frac{1}{2}(q-\sum_{\alpha \in \fd_q}
\lambda(\frac{a}{\alpha^2+\alpha+b}))\\
&= \frac{1}{2}(q-(-K(\lambda;a)-1)). (\text{by Theorem 2 (b)})
\end{align*}
Further, we see that $w(c_4(a)) =
\frac{1}{2}w(c_3(a))$.~~~~~~~~~~~~~~~~~~~~~~~~~~~~~~~~~~~~~~~~~~~~
$\square$
\end{pf*}

Denote for the moment $v_i$ in (11)-(14) by
$v_i=(g_1,g_2,\ldots,g_{N_i})$, for $i=1,2,3,4.$ Fix $i(i=1,2,3,4)$,
and let $u=(u_1,\ldots,u_{N_i}) \in \fd_2^{N_i}$, with $\nu_\beta$
1's in the coordinate places where $g_l = \beta$, for each $\beta
\in \fd_q$. Then we see from the definition of the code $C_i$(cf.
(19)) that $u$ is a codeword with weight $j$ if and only if $\dis
\sum_{\beta \in \fd_q} \nu_\beta = j$ and $\dis \sum_{\beta \in
\fd_q} \nu_\beta \beta = 0$(an identity in $\fd_q$). As there are
$\prod_{\beta \in \fd_q} {n_i(\beta) \choose \nu_\beta}$(cf. Remark
3) many such codewords with weight $j$, we obtain the following
result.

\begin{proposition}
Let  $\{C_{i,j}\}_{j=0}^{N_i}$ be the weight distribution of
$C_i$, for each $i=1,2,3,4,$ where $C_{i,j}$ denotes the frequency
of the codewords with weight $j$ in $C_i$. Then
\begin{equation}
C_{i,j} = \sum \prod_{\beta \in \fd_q} {n_i(\beta) \choose
\nu_\beta},
\end{equation}
where the sum runs over all the sets of integers
$\{\nu_\beta\}_{\beta \in \fd_q}$ $(0 \leq \nu_\beta \leq
n_i(\beta))$, satisfying
\begin{equation}
\sum_{\beta \in \fd_q} \nu_\beta = j, ~and ~~\sum_{\beta \in
\fd_q} \nu_\beta\beta = 0.
\end{equation}
\end{proposition}

\begin{corollary}
Let  $\{C_{i,j}\}_{j=0}^{N_i}$ be the weight distribution of
$C_i$, for $i=1,3.$  Then we have $C_{i,j} = C_{i,N_i-j}$, for all
$j$, with $0 \leq j \leq N_i$.
\end{corollary}
\begin{pf*}{Proof.}
Under the replacements $\nu_\beta \rightarrow n_i(\beta)-\nu_\beta$,
for each $\beta \in \fd_q$, the first equation in (30) is changed to
$N_i-j$, while the second one in (30) and the summands in (29) are
left unchanged. Here the second sum in (30) is left unchanged, since
$\dis \sum_{\beta \in \fd_q} n_i(\beta) \beta = 0$, for $i=1,3$~(cf.
(15), (17)). ~~~~~~~~~~~~~~~~~~$\square$
\end{pf*}

Now, we get the formulas in (2), (5), (7), and (9), by applying the
formula in (29) to each $C_i$, using the values of $n_i(\beta)$ in
(15)-(18).

From now on, we will assume that $r \geq 3$, for $i=1,2$, and hence,
for $i=1,2,3,4$, every codeword in $C_i^\bot$ can be written as
$c_i(a)$, for a unique $a \in \fd_q$ (cf. Theorem 5). Then we apply
the Pless power moment identity in (24) to each $C_i^\bot$, for
$i=1,2,3,4,$ in order to obtain the results in Theorem 1(cf. (1),
(2), (4)-(9)) about recursive formulas. Then the left hand side of
that identity in (24) is equal to
\begin{equation}
\sum_{a \in \fd_q^*} w(c_i(a))^h,
\end{equation}
with the $w(c_i(a))$ in each case given by (25)-(28).

\end{document}